\documentclass[a4paper,12pt]{article} 

\usepackage{amsmath,amssymb}
\usepackage[usenames]{color}
\usepackage{graphicx}
\usepackage{amsfonts}
\usepackage{comment}
\usepackage{pxrubrica}

\newtheorem{thm}{Theorem}
\newtheorem{prop}[thm]{Proposition}
\newtheorem{lmm}[thm]{Lemma}
\newtheorem{defn}{Definition}

\newtheorem{rmk}{Remark}[section]

\numberwithin{equation}{section}

\newcommand{\qed}{\hfill\ensuremath{\square}\\}%
\newcommand{\hf}{\frac{1}{2}}

\newcommand{\pr}{\par \vspace{0mm} \noindent {\bf [Proof]} \quad}
\newcommand{\prend}{\hfill \qed}

\newcommand{\Vir}{{\it Vir}}
 
\newcommand{\1}{{\bf 1}} 

\newcommand{\wt}{{\rm wt}}
\newcommand{\Aut}{{\rm Aut}}

\newcommand{\Rad}{{\rm rad}}

\newcommand{\bC}{{\mathbb C}} 
 
\newcommand{\bZ}{{\mathbb Z}} 
\newcommand{\bN}{{\mathbb N}}

\newcommand{\bM}{{\mathbb M}}

\definecolor{skyblue}{rgb}{0.5,0.5,1}
\definecolor{fadegreen}{rgb}{0.9,1,0.9}

\begin{document}
\title{Borcherds's Lie algebra and $C_2$-cofiniteness of vertex operator algebras of moonshine type}
\author{
\begin{tabular}{c}
Masahiko Miyamoto
\footnote{Partially supported
by the Grants-in-Aids
for Scientific Research, No.21K03195 and No.18K18708, The Ministry of Education,
Science and Culture in Japan
and Academia Sinica in Taiwan} \cr
Institute of Mathematics, University of Tsukuba
\end{tabular}}
\date{}
\maketitle

\begin{abstract}
We precisely determined an $\bN$-graded structure of Zhu's poisson algebra $V/C_2(V)$ of vertex operator algebras $V$ of moonshine type.
\end{abstract}

\section{Introduction}
In these two decades, there was good progress in classifying holomorphic vertex operator algebras (shortly VOA) of CFT-type with a central charge of $24$. Such VOAs were classified under the assumption of strongly regular VOA of CFT-type except for the moonshine type. One of the remaining problems is to prove the uniqueness of VOAs of moonshine type. \\

\noindent
{\bf Uniqueness conjecture of moonshine type in \cite{FLM}} \\
{\it Let $V$ be a vertex operator algebra with a nonsingular invariant bilinear form
$\langle , \rangle$ and a central charge $24$.
If its character $q^{-1}\sum_{n\in \bZ}^{\infty} \dim V_n q^n$ is $j(\tau)-744=q^{-1}+0+196884q+\cdots$,
then $V$ is isomorphic to $V^{\natural}$. }\\

With additional assumptions, several papers have proved the uniqueness, e.g., \cite{DGH}, \cite{DGL}, \cite{T}, \cite{LY}. For example, if there is an isomorphism of Zhu's poisson algebras $\phi: V/C_2(V)\to V^{\natural}/C_2(V^{\natural})$ preserving
inner products $\langle, \rangle$, and grades, then it is not difficult to
prove $V\cong V^{\natural}$ as we will explain in Appendix.
On the other hand, what can we get from the only assumptions in the uniqueness conjecture?
This paper aims to show the precise graded structure of $V/C_2(V)$ of such a VOA $V$.
From now on, let $V$ be a VOA satisfying the conditions in the uniqueness conjecture.
One of the known significant results is that Borcherds's Lie algebra
$B(V)$ of $V$ is isomorphic to Monster Lie algebra $B(V^{\natural})$, which was proved by Borcherds \cite{B3}.
Borcherds has also proved that $B(V)$ is a generalized Kac-Moody algebra
(or Borcherds-Kac-Moody algebra) whose simple roots are $\{(1,-1),(1,1),(1,2),(1,3),... \}\subseteq {\rm II}_{1,1}$ with multiplicities $\dim V_{m+1}$ for each simple root $(1,m)$.
The critical point of his proof is that he has used only
the nonsingular invariant bilinear form
$\langle , \rangle$ of $V^{\natural}$, the fact that the character of $V^{\natural}$ is $J(\tau)$, the Borcherds identity on $J(\tau)$,
and the central charge of $V^{\natural}$ is $24$, see \cite{B} and \cite{J}.
Therefore, it is not difficult to show $B(V)\cong B(V^{\natural})$ as they mentioned 
in \cite{B} and \cite{J}.

As an application of this isomorphism, we will show the following theorem.

\begin{thm}
$C_2(V)=\sum_{n=5}^{\infty}V_n+L(-1)V$. 
In particular, V is $C_2$-cofinite.
\end{thm}

The $C_2$-cofiniteness of $V^{\natural}$ was initially proved by Dong, Griess, and Hohn \cite{DGH}. They showed that $V^{\natural}$ contains a framed structure $L(\hf,0)^{\otimes 48}$,
where $L(\hf,0)$ is one of discrete series with a central charge $1/2$, and it was proved to be $C_2$-cofinite by direct calculation (see \cite{DLM}). \\

\noindent
{\bf Acknowledgement} \qquad
The author thanks T.~Abe, C.~H.~Lam, and H.~Yamauchi for their suggestions and valuable comments.
He also thanks Academia Sinica in Taiwan for the hospitality.

\section{Borcherds's Lie algebra of VOA of moonshine type}
Set
$\tilde{V}=V\otimes V_{{\rm II}_{1,1}}$, where ${\rm II}_{1,1}$ denotes an even unimodular Lorentian lattice of rank $2$ with Gramm matrix $\begin{pmatrix}0, -1\cr -1,0
\end{pmatrix}$
and $V_{{\rm II}_{1,1}}$ denotes a lattice vertex algebra defined by
${\rm II}_{1,1}$.
Then $\tilde{V}$ has a ${\rm II}_{1,1}$-graded structure:
$$V\otimes V_{{\rm II}_{1,1}}=\oplus_{(m,n)\in {\rm II}_{1,1}} V\otimes V_{{\rm II}_{1,1}}^{(m,n)}=\oplus_{(m,n)\in {\rm II}_{1,1}}V\otimes M(1)^{\otimes 2}e^{(m,n)},$$
where $M(1)$ denotes a Heisenberg VOA of rank one.
Let $\1$ and $\1'$ denote Vaccumes of $V$ and $V_{{\rm II}_{1,1}}$, respectively. Then $\tilde{\1}=\1\otimes \1'$ is a Vaccume of $\widetilde{V}$. Also
$\{L(n)\mid n\in \bZ\}$, $\{L'(n)\mid n\in \bZ\}$, and $\{\tilde{L}(n)=L(n)+L'(n)\mid n\in \bZ\}$ denote
Virasoro operators of $V$, $V_{{\rm II}_{1,1}}$, and $\tilde{V}$, respectively.
From now on, $U(\Vir)$ denotes the universal envelopping algebra of Virasoro
algebra $\Vir=\oplus_{n\in \bZ} \bC L(n)+\bC$ with central charge $24$.

Set
$$P^k=\{u\in \tilde{V} \mid \tilde{L}(0)u=ku, \tilde{L}(n)u=0 \,\, \forall n> 0\}.$$
$P^1$ is called the space of physical states.
From the property of $0$-product in $\widetilde{V}$,
$P^1/\tilde{L}(-1)P^0$ becomes a Lie algebra with an invariant bilinear form induced from the invariant bilinear form
$\langle \cdot, \cdot\rangle_{\tilde{V}}$ on $\tilde{V}$.
Then we define a Borcherds's Lie algebra $B(V)$ of a VOA $V$
as the quotient of $P^1$ by the radical (or the null space) ${\rm rad}(P^1)$ of $\langle\cdot,\cdot\rangle_{\tilde{V}}$.
We note $\tilde{L}(-1)P^0\subseteq {\rm rad}(P^1)$.
$P^1_{(m,n)}$ denotes the subspace of $P^1$ of degree
$(m,n)\in {\rm II}_{1,1}$ and set
$B(V)^{(m,n)}=P^1_{(m,n)}/\Rad(P^1_{(m,n)})$.

We recall the following result from \cite{B3}, which is originally mentioned
by Goddard and Thorn and also called the Goddard-Thorn theorem \cite{GT}.

\begin{thm}[The no-ghost theorem]
Suppose that $V=\oplus_{n=0}^{\infty} V_n$ is a VOA of central charge 24
with a nonsingular invaariant bilinear form $\langle ,\rangle$ and $G\leq \Aut(V)$.
$V\otimes V_{{\rm II}_{1,1}}$ inherits an action of $G$ from the action of $G$ on $V$ and the trivial action of $G$ on $V_{{\rm II}_{1,1}}$ and $\bC^2$.
Then the quotient $B(V)^{(m,n)}$ of $P^1_{(m,n)}$ by the nullspace $\Rad(P^1_{(m,n)})$ of $\langle , \rangle$ is naturally isomorphic,
as a $G$-module with an invariant bilinear form, to $V_{1+mn}$ if $mn\not=0$ and to $V_1\oplus \bC^2$ if $(m,n)=(0,0)$ and
$B(V)^{(m,n)}=0$ for else.
\end{thm}

In the case of the moonshine VOA $V^{\natural}$, since it has a huge automorphism group
$\bM$ called the monster group, an isomorphism $B(V^{\natural})^{(m,n)}\cong V_{mn+1}$ as $\bM$-modules gives strong information. However, in our setting, 
we cannot use a nontrivial automorphism group for a general VOA $V$ satisfying the conditions in Theorem 1.
Therefore, the no-ghost theorem only says
$$\dim B(V)^{(m,n)}=\dim V_{mn+1}.$$
In order to use the no-ghost theorem more effectively, we will introduce a large group, which
may not be an automorphism group of $V$, but an orthogonal linear transformation group of the inner product space $(V,\langle \cdot,\cdot\rangle)$.
Since $P^1$ and $B(V)$ are defined by Virasoro operators
$\{\tilde{L}(n)=L(n)+L'(n):n\in \bZ\}$, we will introduce a group
$$G=\{ g=\prod_m g_m\in \prod_m O(V_m) \mid g_mL(n)=L(n)g_{m+n}
\, \forall n,m\in \bZ\},$$
which consists of orthogonal linear transformations of an inner product space $V$, preserving the grades and
commuting with the actions of Virasoro operators $\{L(n)\mid n\in \bZ\}$ of $V$.
We extend it to an orthogonal linear treansformation group $G\otimes 1_{V_{{\rm II}_{1,1}}}$ of the inner product space $(\tilde{V}, \langle,\rangle_{\tilde{V}})$.
Set
$$V_k^p=\{v\in V_k \mid L(n)v=0 \,\, {}^\forall n\geq 1\},$$
which is the set of $\Vir$-primary states in $V$ of conformal weight $k$.
Then since $\langle, \rangle$ has the nonsingular invariant bilinear form,
we have $\langle U(\Vir)V_m^p, U(\Vir)V_n^p\rangle=0$ for $m\not=n$ and
$V=\oplus_{n\in \bZ} U(\Vir)V_n^p$. Every $V_m^p$ is a simple $G$-module since any two nonzero elements $v, w$ of $V_m^p$
have the same conformal weight for $\Vir$ and
there is $u\in V$ such that $o(u)v=w$, where $o(u)$ denotes the grade preserving operator of $u$. On the other hand, if $g\in O(V_k^p)$, then since the central charge of
$V$ is 24, the subspace $(U(\Vir)V_k^p)_m$ of weight $m$ is a direct sum of copies of $V_k^p$ and so we can extend $g$ to $O(U(\Vir)V_k^p)$ by $g(L(-m)u)=L(-m)g(u)$ by induction.
Therefore we have $G\cong \prod_{m=0}^{\infty}O(V_m^p)$.

From the definition of physical states, $P^1$ and $B(V)$ are invariant under the actions of $G$, that is, they are $G$-modules and
$P^1_{(m,n)}$ and $B(V)^{(m,n)}$ are $G$-submodules of them.
From now on, we will consider a projection
$$\pi_k:V\otimes V_{{\rm II}_{1,1}} \to U(\Vir)V_k^p\otimes V_{{\rm II}_{1,1}}$$
for each $k\in \bZ$ and
set $P^1_{(m,n),k}:=\pi_k(P^1_{(m,n)})$ and $B(V)^{(m,n),k}:=\pi_k(B(V)^{(m,n)})$.
We note $\pi_k(V)=U(\Vir)V_k^p$.

As one of the steps of the proof of Theorem 1, we will prove
$V_5\subseteq C_2(V)$ and $V_4\subseteq L(-1)V_3+(V_2)_0(V_3)$ in the next section.
Since $L(-1)V, L(-3)V, L(-4)V,... \subseteq C_2(V)$, we have
$V_4\subseteq V_4^p+L(-2)V_2^p+L(-2)^2V_0+C_2(V)$, $V_5\subseteq L(-2)V_3^p+V_5^p+C_2(V)$, and
$V_6\subseteq V_6^p+L(-2)V_4+C_2(V)$.
Therefore we first study the structures of $P^1_{(2,2),3}+P^1_{(2,2),5}$ for $L(-2)V_3^p+V_5$ and those of
$P^1_{(2,3)}$ for $L(-2)^3V_0+L(-2)^2V_2^p+L(-2)V_4^p+V_6^p$ in this section.

We need the following results.

\begin{lmm}\label{Afewcase}
For $m>0$ and $m+n\leq 5$, $B(V)^{(m,n)}\cong V_{mn+1}$ or $0$ as $G$-modules.
\end{lmm}

\pr By direct calculation, since $\wt(e^{(m,n)})=-mn$, it is easy to see
$P^1_{(m,n),mn+1}=V_{mn+1}^p\otimes e^{(m,n)}$ and $\Rad(P^1_{(m,n),mn+1})=0$.
We also have
$P^1_{(m,n),mn}=V_{mn}^p\otimes \delta^{(m,-n)}_{(-1)}e^{(m,n)}+
\tilde{L}(-1)(V_{mn}^p\otimes e^{(m,n)})$ and
$\Rad(P^1_{(m,n),mn})=\tilde{L}(-1)(V_{mn}^p\otimes e^{(m,n)})$,
where $\delta^{(m,n)}_{(-1)}e^{(a,b)}$ denotes $(m,n)_{(-1)}e^{(a,b)}$
for $(m,n)\in {\rm II}_{1,1}$ by viewing $(m,n)$ as an element in $(M(1)^{\otimes 2})_1
\subseteq V_{II_{1,1}}$. 
Therefore, we have $B(V)^{(m,n),mn+1}\cong V_{mn+1}^p$ and $B(V)^{(m,n),mn}\cong V_{mn}^p$ as $G$-modules.
Furthermore, we have obtained $\dim B(V)^{(m,n)}=\dim V_{mn+1}$ from the 
no-ghost theorem.
Using the above results, we get $B(V)^{(1,-1)}\cong V_0$,
$B(V)^{(1,1)}\cong V_0\oplus V_2^p\cong V_2$, and
$B(V)^{(1,2)}\cong V_0\oplus V_2^p\oplus V_3^p\cong V_3$.
Since $\dim B(V)^{(1,3),0}\leq
\dim (U(\Vir)\1\otimes M(1)^{\otimes 2}e^{(1,3)})_1 <\!\!<\dim V_2^p=196883$, we
also obtain 
$B(V)^{(1,3)} \cong V_4^p+V_3^p+2V_2^p+2V_0\cong V_4$ as $G$-modules.
For $B(V)^{(1,4)}$, $B(V)^{(2,2)}$, and $B(V)^{(2,3)}$, we will prove them 
after Proposition 10.
\prend

Although we will not use it, the following remark comes from Borcherds's proof of the no-ghost theorem since he used only elements in $\Vir \otimes M(1)^{\otimes 2}$ and the actions of these elements preserve $V_k^p$ in his proof.

\begin{rmk}[Borcherds's proof in \cite{B3}]
$(B(V)^{(m,n)}, \langle,\rangle) \cong (V_{mn+1},\langle,\rangle)$ as $G$-modules.
\end{rmk}

Borcherds has also shown in \cite{B3} that Borcherds's Lie algebra
$B(V)=P^1/{\rm rad}(P^1)$ is a generalized Kac-Moody algebra.

\begin{thm}[\cite{B3}]
If $V$ satisfies the conditions of the uniqueness conjecture, then the simple roots of Lie algebra $B(V)$ are the vectors
$(1,m)$ $(m=-1$ or $m>0)$ with multiplicity $\dim V_{m+1}$.
\end{thm}

The most useful result is that
$\{B(V)^{(1,n)}: n=-1,1,2,...\}$ generates all positive root spaces
$\oplus_{m=1}^{\infty}\oplus_{n\in \bZ} B(V)^{(m,n)}$. In particular, we will use
$$B(V)^{(2,m)}=[B(V)^{(1,m+1)}, B(V)^{(1,-1)}]
+\sum_{k=1}^{m-1}[B(V)^{(1,k)}, B(V)^{(1,m-k)}] \mbox{ for all $m$}. \eqno{(A)}$$
In other words, coming back to physical states, we have
$$P^1_{(2,m)} =(P^1_{(1,-1)})_0(P^1_{(1,m+1)})
+\sum_{k=1}^{m-1}(P^1_{(1,k)})_0(P^1_{(1,m-k)})+{\Rad}(P^1_{(2,m)}). \eqno{(B)}$$

Now, we can explain one of our arguments.

\begin{prop}\label{C2}
If $n\geq 5$ and $L(-2)V_{n-2}\subseteq C_2(V)$,
then $V_n\subseteq C_2(V)$.
\end{prop}

\pr
Assume that $n=2m+1\geq 5$ is odd.
If $v\in V_{2m+1}^p$, then as we have explained,
$v\otimes e^{(2,m)}\in P^1_{(2,m)}\backslash \Rad(P^1)$.
By (B), there are $\sum_a u^{k,a,i}\otimes \alpha^{k,a,i}e^{(1,k)}\in P^1_{(1,k)}, \sum_b w^{m-k,b,j}\otimes \beta^{m-k,b,j}e^{(1,m-k)}\in P^1_{(1,m-k)}$,
$\sum_i w^{m+1,i}\otimes \beta^{m+1,i}e^{(1,m+1)}\in P^1_{(1,m+1)}$, and
$\sum_i z^i\otimes \delta^i e^{(1,2m)}\in \Rad(P^1_{(2,m)})$ with
$u^{k,a,i},w^{m-k,b,i},z^i\in V$ and $\alpha^{k,i}e^{(1,a)},\beta^{m-k,i}e^{(1,a)}, \delta^ie^{(1,2m)}\in V_{{\rm II}_{1,1}}$ such that
$$\begin{array}{rl}
v\otimes e^{(2,m)}=&(1\otimes e^{(1,-1)})_0 (\sum_i w^{m+1,i}\otimes \beta^{m+1,i}e^{(1,m+1)})\cr
&+\sum_{k=1}^{m-1}\sum_j \sum_a (u^{k,a,j}\otimes \alpha^{k,a,j}e^{(1,k)})_0(\sum_b w^{m-k,b,j}\otimes \beta^{m-k,b,j}e^{(1,m-k)})\cr
&+\sum_j z^j\otimes \delta^je^{(1,2m)}.
\end{array}$$
In particular, $v$ is a linear sum of $\{w^{m+1,i}, u^{k,a,i}_h w^{m-k,b,s}, z^j\mid k,i,j,a,b,s,h \}$.
Since $\sum_a u^{k,a,i}\otimes \alpha^{k,a,i}e^{(1,k)}\in P^1_{(1,k)}$ and $\sum_b w^{m-k,b,i}\otimes \beta^{m-k,b,i}e^{(1,m-k)}\in P^1_{(1,m-k)}$,
we have $\wt(u^{k,a,i})\leq k+1$ and
$\wt(w^{m-k,b,i})\leq m-k+1$. Hence we have
$\wt(u^{k,i})+\wt(w^{m-k,i})\leq m+2$.
Also since $\Rad(P^1_{(2,m),2m+1})=0$, we may choose
$z^i\in (V_{2m+1}^p)^{\perp}=U(\Vir)\sum_{k=0}^{2m}V_{k}$. We hence have $\wt(z^i)\leq 2m+1$.
Furthermore, since $\wt(v)=2m+1$, in order to get an element $v$
of $V_{2m+1}^p$ as a linear sum
of $\{u^{k,a,i}_h w^{m-k,b,j},z^i\mid k,i,j,a,b,h\}$, we have $h\leq -m\leq -2$ and
$v\in \pi_{2m+1}(C_2(V))$ for every $v\in V_{2m+1}^p$.
Therefore, $V_{2m+1}^p\subseteq \pi_{2m+1}(C_2(V))$ and $V_{2m+1}^p\subseteq
C_2(V)+((V_{2m+1})^{\perp})_{2m+1}\subseteq C_2(V)+L(-2)V_{2m-1}$.
We hence have
$V_{2m+1}\subseteq C_2(V)+L(-2)V_{2m-1} \subseteq C_2(V)$ by the assumption.

We next assume $n=2m\geq 6$. By (B), we have: 
$$B(V)^{(2,m)}=[B(V)^{(1,-1)},B(V)^{(1,m+1)}]+\sum_{k=1}^{m-1}[B(V)^{(1,k)},B(V)^{(1,m-k)}].$$
Therefore, for $v\in V^p_{2m}$, since $v\otimes \delta^{(2,-m)}_{(-1)}e^{(2,m)}\in P^1$, we have: 
$$\begin{array}{l}
v\otimes \delta^{(2,-m)}_{(-1)} e^{(2,m)}\in
(\bC\1\otimes e^{(1,-1)})_0((\oplus_{n=0}^{m+2}V_n)\otimes M(1)^{\otimes 2}e^{(1,m+1)})\cr
\mbox{}\qquad \qquad + \sum_{k=1}^{m-1}((\oplus_{h=0}^{k+1}V_h)\otimes M(1)^{\otimes 2}e^{(1,k)})_0((\oplus_{h=0}^{m-k+1}V_h)\otimes M(1)^{\otimes 2}e^{(1,m-k)})+{\rm rad}(P^1).
\end{array}$$
Since $||v\otimes \delta^{(2,-m)}_{(-1)}e^{(2,m)}||\not=0$ and
$\{v\otimes \delta^{(2,-m)}_{(-1)} e^{(2,m)}\mid v\in V_{2m}^p \}+\Rad(P^1_{(2,m),2m})=P^1_{(2,m),2m}$,
we can choose an orthogonal basis of $P^1_{(2,m)}$ containing
$v\otimes \delta^{(2,-m)}_{(-1)}e^{(2,m)}$ as a base.
When we express elements in $P^1_{(2,m)}$ as a linear combination of the orthogonal basis, any elements $\xi\in \Rad(P^1_{(2,m)})$ have no coefficients at
$v\otimes \delta^{(2,-m)}_{(-1)}e^{(2,m)}$, since $\langle \Rad(P^1), v\otimes \delta^{(2,-m)}_{(-1)}e^{(2,m)}\rangle=0$.
Therefore, $v\in V^p_{2m}$ is a linear sum of images $\pi_{2m}(u^i_hw^i)$ of the products $u^i_{h_i}w^i$ with $u^i\in \oplus_{h=0}^{k+1}V_h$ and
$w^i\in \oplus_{h=0}^{m-k+1}V_h$ modulo $(V^p_{2m})^{\perp}\cap V_{2m}$,
that is, $v=\pi_{2m}(\sum_i u^i_{h_i}w^i)$.
Since $\wt(u^i)+\wt(w^i)\leq m+2$, we have $h_i\leq -2$ and $v\in \pi_{2m}(C_2(V))$ for every $v\in V_{2m}^p$, which implies
$V_{2m}^p\subseteq C_2(V)+(V_{2m}^p)^{\perp}$.
Since $(V^p_{2m})^{\perp}\cap V_{2m}\subseteq L(-1)V+L(-2)V_{2m-2}$,
we have $V_{2m}\subseteq C_2(V)+L(-2)V_{2m-2}\subseteq C_2(V)$ by the assumption.
\prend

\section{Structure of physical states}
We recall vertex operators $Y(e^{\gamma}, z)$ of $e^\gamma$ for a lattice VOA $V_{{\rm II}_{1,1}}$ and $\gamma\in II_{1,1}$ from \cite{FLM}. Namely,
$$\begin{array}{rl}
Y(e^\gamma,z)=&E^-(-\gamma,z)E^+(-\gamma,z)e^{\gamma}z^{\wt(e^{\gamma})+\gamma}\cr
=&
\exp(\sum_{n=1}^{\infty} \frac{\gamma(-n)}{n}z^n)\exp(\sum_{n=1}^{\infty}-
\frac{\gamma(n)}{n} z^{-n})e^{\gamma}z^{\wt(e^{\gamma})+\gamma}.
\end{array}\eqno{(OP)}$$
For $\mu\in (V_{{\rm II}_{1,1}})_m$, we use the notation
$$Y(\mu,z)=\sum_{k\in \bZ} o_k(\mu) z^{k-m}, \mbox{ where } [L(0)',o_k(\mu)]=ko_k(\mu). $$

From now on, $\overline{A}$ denotes $(A+\Rad(P^1))/\Rad(P^1)$
for a subset $A$ of $P^1$. 

Since $G$ does not preserve products in $V$,
$(P^1_{(1,m),k})_0(P^1_{(1,n),h})$ may not be a $G$-module.
To calculate $(P^1_{(1,m),k})_0(P^1_{(1,n),h})$, let us put $G$-modules
aside.
Let $Z$ be a $G$-submodule of a $G$-module $Y$ and
$X$ a subspace of $Y$.
If $(X+Z)/Z$ contains a $G$-submodule $T$ of $Y/Z$, then we will say that
$X$ covers a $G$-module $T$ modulo $Z$.

\subsection{$G$-homomorphisms and $0$-products by $(P^1_{(1,-1)})_0$ and $(\pi_0(P^1))_0$}
Clearly, $P^1_{(1,-1)}=\bC \1\otimes e^{(1,-1)}$ and its $0$-product
$(\1\otimes e^{(1,-1)})_0=1\otimes (e^{(1,-1)})_0$ is a $G$-homomorphism.
More generally, for $\alpha\in \pi_0(P^1)\subseteq U(\Vir)\1\otimes V_{{\rm II}_{1,1}}$,
its $0$-product $\alpha_0$ is given by Virasoro opertors and operators of elements in $V_{{\rm II}_{1,1}}$, which commute with $G$.
Therefore, $0$-product $\alpha_0$ is also a $G$-homomoprhism.
Hence $[B(V)^{(1,-1)}, B(V)^{(1,k)}]$ and $[\pi_0(B(V)^{(1,m)}), B(V)^{(1,k)}]$ are holomorphic images of $B(V)^{(1,k)}$
as $G$-modules for $m=1,2$ since $\dim \pi_0(V_2)=\dim \pi_0(V_3)=1$.

\subsection{$[B(V)^{(1,1),2}, B(V)^{(1,1),2}]$ and $[B(V)^{(1,1).2}, B(V)^{(1,2),2}+B(V)^{(1,2),3}]$}
Set $X=(P^1_{(1,1),2})_0(P^1_{(1,1),2})$ and $Y=(P^1_{(1,1),2})_0(P^1_{(1,2),2})
+(P^1_{(1,1),2})_0(P^1_{(1,2),3})$. Since \\
$B(V)^{(2,2)}=[B(V)^{(1,-1)},B(V)^{(1,3)}]+[B(V)^{(1,1),0}, B(V)^{(1,1)}]
+\overline{X}$, $X$ covers \\
$B(V)^{(2,2)}/([B(V)^{(1,-1)},B(V)^{(1,3)}]+[B(V)^{(1,1),0}, B(V)^{(1,1)}])$ modulo $\Rad(P^1)$. Similarly,
$Y$ covers
$B(V)^{2,3)}/([B(V)^{(1,-1)},B(V)^{(1,3)}]+[B(V)^{(1,1),0}, B(V)^{(1,2)}]+[B(V)^{(1,1)}, B(V)^{(1,2),0}])$ modulo $\Rad(P^1)$.
By expressing $G$-modules by sums of copies of $V_k^p$ after we prove
 $B(V)^{(m,n)}\cong V_{mn+1}$ for $(m,n)=(2,2),(2,3)$,
we will obtain:

\begin{prop}\label{Cover2} \mbox{}\\
(1) $(V_2^p\otimes e^{(1,1)})_0(V_2^p\otimes e^{(1,1)})$ covers $V_3^p\oplus V_5^p$ modulo $\Rad(P^1)$ and \\
(2) $(V_2^p\otimes e^{(1,1)})_0(V_2^p\otimes \delta^{(1,-2)}_{(-1)}e^{(1,2)})
+(V_2^p\otimes e^{(1,1)})_0(V_3^p\otimes e^{(1,2)})$ covers
$V_0\oplus 2V_2^p\oplus 2V_4^p\oplus V_6^p$ \\
\mbox{}\qquad modulo $\Rad(P^1)$.
\end{prop}

\subsection{$P^1_{(1,1)}$ and $P^1_{(1,2)}$ and their products}
Let us study $(P^1_{(1,1),2})_0P^1_{(1,1),2}$
and
$(P^1_{(1,1),2})_0P^1_{(1,2),3}$ elementwisely.

For $v,u\in V_2^p$ and $w\in V_3^p$, we define
$$\begin{array}{ll}
\Phi^1(v,u)=(v\otimes e^{(1,1)})_0(u\otimes e^{(1,1)})&\in P^1_{(2,2)} \cr
\Phi^2(v,u)=(v\otimes e^{(1,1)})_0(u\otimes \delta^{(1,-2)}_{(-1)}e^{(1,2)})
&\in P^1_{(2,3)} \qquad\mbox{ and}\cr
\Phi^3(v,w)=(v\otimes e^{(1,1)})_0(w\otimes e^{(1,2)})&\in P^1_{(2,3)}.
\end{array}$$

Using $Y(v\otimes e^{(1,1)},z)=Y(v,z)\otimes \sum_{k\in \bZ} o_{k}(e^{(1,1)})z^{k+1}$
and $o_k(e^{(1,1)})(V_{{\rm II}_{1,1}})_m\subseteq (V_{{\rm II}_{1,1}})_{m+k}$,
we have expressions:
$$\begin{array}{rl}
\Phi^1(v,u)=&v_{-2}u\otimes o_0(e^{(1,1)})e^{(2,2)}
+v_{-1}u\otimes o_1(e^{(1,1)})e^{(2,2)}+v_0u\otimes
o_2(e^{(1,1)})e^{(2,2)}+\cdots \cr
\Phi^2(v,u)=&v_{-4}u\otimes o_{-1}(e^{(1,1)})\delta^{(1,-2)}_{(-1)}e^{(2,3)}
+v_{-3}u\otimes o_{0}(e^{(1,1)})\delta^{(1,-2)}_{(-1)}e^{(2,3)}\cr
&+v_{-2}u\otimes o_1(e^{(1,1)})\delta^{(1,-2)}_{(-1)}e^{(2,3)}
+v_{-1}u\otimes o_2(e^{(1,1)})\delta^{(1,-2)}_{(-1)}e^{(2,3)}\cr
&+v_0u\otimes o_3(e^{(1,1)})\delta^{(1,-2)}_{(-1)}e^{(2,3)}+\cdots \cr
\Phi^3(v,w)=&v_{-3}u\otimes o_0(e^{(1,1)})e^{(2,3)}+
v_{-2}u\otimes o_1(e^{(1,1)})e^{(2,3)}
+v_{-1}u\otimes o_2(e^{(1,1)})e^{(2,3)}\cr
&+v_0u\otimes
o_3(e^{(1,1)})e^{(2,3)}+\cdots
\end{array} \eqno{(C)}$$
We note that
$\{ o_k(e^{(1,1)})e^{(2,2)}, \,
o_k(e^{(1,1)})\delta^{(1,-2)}_{(-1)}e^{(2,3)}, \,
o_k(e^{(1,1)})e^{(2,3)}\}\,\, (\subseteq V_{{\rm II}_{1,1}})$
do not depend on the choices of $v,u,w$.

We next study their images by the projections
$\pi_k: \tilde{V} \to U(\Vir)V_k^p\otimes V_{{\rm II}_{1,1}}$.
Clearly, $v_{3-k}u, v_{4-k}w\in V_k$ for $v,u\in V_2^p, w\in V_3^p$. 
Since $\wt(v_mu)=3-m$ and $\wt(v_mw)=4-m$, 
$\pi_{k}(v_{m}u)=0$ for $m>3-k$ and
$\pi_{k}(v_mw)=0$ for $m>4-k$.
We note that if $\pi_k(v_{3-k}u)=0$, then $\pi_k(v_{3-k-j}u)=0$ for all $j\in \bN$,
since all primary states of $\pi_k(V)=U(\Vir)V_k^p$ are in $V_k^p$.
Similarly, if $\pi_k(v_{4-k}w)=0$, then $\pi_k(v_{4-k-j}w)=0$ for all $j\in \bN$.

\begin{lmm}
For each pair $(k,m)\in \bN\times \bN$,
there are operators $Q^{2,2}_k(-m), Q^{2,3}_k(-m)\in U(\Vir)$ of degree $-m$ such that
$\pi_{k}(v_{3-k-m}u)=Q^{2,2}_k(-m)\pi_{k}(v_{3-k}u)$ and
$\pi_{k}(v_{4-k-m}w)=Q^{2,3}_k(-m)\pi_{k}(v_{4-k}w)$
for any $v,u\in V_2^p$ and $w\in V_3^p$.
\end{lmm}

\pr
Since all primary states of $\pi_{k}(V)$ are in $V_k^p$,
if $\pi_k(v_{3-k-m}u)\not=0$ for some $m$, then $\pi_k(v_{3-k}u)\not=0$ and $\pi_k(v_{3-k-m}u)\in U(\Vir)\pi_k(v_{3-k}u)$.
Therefore there is $Q^{2,2}_k(-m)\in U(\Vir)$ of degree $m$
such that $\pi_k(v_{3-k-m}u)=Q^{2,2}_k(-m)\pi_k(v_{3-k}u)$.
Since the central charge is $24$, we can express $U(\Vir)V_k^p$ as a linear sum of $L(-n_1)...L(-n_k)V_k^p$ with $n_1\geq ... \geq n_k\geq 1$
uniquely.
Since $[L(k), v_{-s}]=(s+k)v_{-s+k}$ holds for any $v\in V_k^p$ and $s,k\in \bZ$, 
we can choose $Q^{2,2}_k(-m)$ so that $Q^{2,2}_k(-m)$
does not depend on the choice of $v,u\in V_2^p$. Similarly, if
$\pi_k(v_{4-k-m}w)\not=0$ for some $m$, then $\pi_k(v_{4-k}w)\not=0$ and
there is $Q^{2,3}_k(-m)\in U(\Vir)$ of degree $m$ such that
$\pi_k(v_{4-k^m}w)=Q^{2,3}_k(-m)\pi_k(v_{4-k}w)$ and $Q^{2,3}_k(-m)$ does not depend on the choice of $v\in V_2^p$ and $w\in V_3^p$.
\prend

We will show the precise formula for several $Q^{a,b}(-m)$, which we need.

\begin{lmm}\label{Direct}
For $v,u\in V_2^p$ and $w\in V_3^p$, set $X=\pi_{0}(v_3u)$, $Y=\pi_2(v_1u)$, and $Z=\pi_3(v_1w)$. Then we have:
$$\begin{array}{l}
\pi_k(v_{2-k}u)=\frac{1}{2}L(-1)\pi_k(v_{3-k}u), \vspace{1mm} \cr
\pi_k(v_{3-k}w)=\frac{k-1}{2k}L(-1)\pi_k(v_{4-k}w), \vspace{1mm} \cr
\pi_{0}(v_1u)=\frac{1}{6}L(-2)X, \vspace{1mm} \cr
\pi_{0}(v_0u)=\frac{1}{12}L(-3)X, \vspace{1mm} \cr
\pi_{0}(v_{-1}u)=\frac{3}{71}L(-4)X+\frac{11}{852}L(-2)^2X, \vspace{1mm} \cr
\pi_{0}(v_{-2}u)=\frac{2}{71}L(-5)X+\frac{11}{852}L(-3)L(-2)X, \vspace{1mm} \cr
\pi_{0}(v_{-3}u)=\frac{1}{196883}\{3492L(-6)X+\frac{15623}{12}L(-4)L(-2)X+\frac{1271}{2}L(-3)^2X
+124L(-2)^3X\}, \vspace{1mm} \cr
\pi_2(v_0w)=\frac{6}{41}L(-2)Y+\frac{1}{164}L(-1)^2Y, \hfill \mbox{ and \qquad \qquad } \vspace{1mm} \cr
\pi_{3}(v_{-2}u)=\frac{17}{141}L(-2)Z+\frac{11}{94}L(-1)^2Z.
\end{array}$$
In particular, we have
$L(-2)Z\in \pi_3(C_2(V))$, $L(-2)^3\1\in \pi_0(C_2(V))$, and
$L(-2)Y\in \pi_2((V_2)_0(V_3)+L(-1)V_3)$.
\end{lmm}

\pr
We will prove all equations by direct calculation. 
As we mentioned, $\pi_k(v_{2-k}u)=aL(-1)\pi_k(v_{3-k}u)$ and $\pi_k(v_{3-k}w)=bL(-1)\pi_k(v_{4-k}w)$ for some $a,b\in \bC$. Then we obtain 
$L(1)\pi_k(v_{2-k}u)=k\pi_k(v_{3-k}u)$,
$L(1)aL(-1)\pi_k(v_{3-k}u)=2ka\pi_k(v_{3-k}u)$, and
$L(1)\pi_k(v_{3-k}w)=(k-1)\pi_k(v_{4-k}w)=2bk \pi_k(v_{4-k}w)$.
Hence we have $a=\frac{1}{2}$ and $b=\frac{k-1}{2k}$.

Since $\pi_0(v_1u)\in L(-2)V_0$, there is $a\in \bC$ such that $\pi_0(v_1u)=aL(-2)\pi_0(v_3u)$ since as we showed, if $\pi_0(v_3u)=0$, then $\pi_0(v_1u)=0$. 
Then $L(2)\pi_0(v_1u)=2\pi_0(v_3u)$ and $L(2)aL(-2)\pi_0(v_3u)=a(4L(0)+12)\pi_0(v_3u)=12a\pi_0(v_3u)$.
Hence we have $a=\frac{1}{6}$ and $\pi_0(v_1u)=\frac{1}{6}L(-2)\pi_0(v_3u)$.
Since $\wt(\pi_0(v_0u))=3$, there is $a\in \bC$ such that
$\pi_0(v_0u)=aL(-3)\pi_0(v_3u)$.
Then $\frac{1}{3}L(-2)\pi_0(v_3u)=2\pi_0(v_1u)=L(1)\pi_0(v_0u)=4aL(-2)\pi_0(v_3u)$
and so $a=\frac{1}{12}$.
For $\pi_0(v_{-1}u)$, there are $a,b\in \bC$ such that
$\pi_0(v_{-1}u)=aL(-4)\pi_0(v_3u)+bL(-2)^2\pi_0(v_3u)$.
Then since $\frac{1}{4}L(-3)\pi_0(v_3u)=3\pi_0(v_0u)=L(1)\pi_0(v_{-1}u)
=5aL(-3)\pi_0(v_3u)+3bL(-3)\pi_0(v_3u)$, we have $5a+3b=\frac{1}{4}$.
Furthermore, $\frac{2}{3}L(-2)\pi_0(v_3u)=4\pi_0(v_1u)=L(2)\pi_0(v_{-1}u)
=6aL(-2)\pi_0(v_3u)+20bL(-2)\pi_0(v_3u)+12bL(-2)\pi_0(v_3u)$ and so
we have $6a+32b=\frac{2}{3}$. Solving the simultaneous equation, we have
$a=\frac{3}{71}$ and $b=\frac{11}{852}$.
For $\pi_0(v_{-2}u)$, there are $a,b\in \bC$ such that $\pi_0(v_{-2}u)=aL(-5)\pi_0(v_3u)+bL(-3)L(-2)\pi_0(v_3u)$.
Then $\frac{12}{71}L(-4)\pi_0(v_3u)+\frac{44}{852}L(-2)^2\pi_0(v_3u)=4\pi_0(v_{-1}u)=L(1)\pi_0(v_{-2}u)=5aL(-4)\pi_0(v_3u)+4bL(-2)^2\pi_0(v_3u)$.
Hence we have $a=\frac{2}{71}$ and $b=\frac{11}{852}$.
For $\pi_0(v_{-3}u)$, there are $a,b,c,d\in \bC$ such that
$\pi_0(v_{-3}u)=aL(-5)\pi_0(v_{-3}u)+bL(-4)L(-2)\pi_0(v_3u)+cL(-3)^2\pi_0(v_3u)+dL(-2)^3\pi_0(v_3u)$.
Then
$$\begin{array}{l}
\frac{10}{71}L(-5)\pi_0(v_3u)+\frac{55}{852}L(-3)L(-2)\pi_0(v_3u)
=5\pi_0(v_{-2}u)=L(1)\pi_0(v_{-3}u)\cr
\mbox{}\qquad =7aL(-5)\pi_0(v_3u)+5bL(-3)L(-2)\pi_0(v_3u)+\{8cL(-3)L(-2)+4cL(-5)\}\pi_0(v_3u)\cr
\mbox{}\qquad +\{6dL(-5) +9dL(-3)L(-2)\}\pi_0(v_3u).
\end{array}$$
Therefore have $7a+4c+6d=\frac{10}{71}$ and $5b+8c+9d=\frac{55}{852}$.
Then $L(2)\pi_0(v_{-3}u)=6\pi_0(v_{-1}u)=\frac{18}{71}L(-4)\pi_0(v_3u)
+\frac{66}{852}L(-2)^2\pi_0(v_3u)$
and $6\pi_0(v_{-1}u)=8aL(-4)\pi_0(v_3u)+\{6bL(-2)^2+12bL(-4)\}\pi_0(v_3u)
+60dL(-2)^2\pi_0(v_3u)$.
Therefore we have $8a+12b+10c=\frac{18}{71}$ and $6b+60d=\frac{66}{852}$.
Solving the simultaneous equation, we have
$a=\frac{3492}{196883}$, $b=\frac{15623}{2362596}$, $c=\frac{1271}{393766}$, and
$d=\frac{124}{196883}$.

For $\pi_2(v_0w)$, there are $a,b\in \bC$ such that $\pi_2(v_0w)=aL(-2)\pi_2(v_2w)
+bL(-1)^2\pi_2(v_2w)$. Then $\frac{1}{2}L(-1)\pi_2(v_2w)=2\pi_2(v_1w)=L(1)\pi_2(v_0w)
=3aL(-1)\pi_2(v_2w)+10bL(-1)\pi_2(v_2w)$ and so $3a+10b=\frac{1}{2}$.
Furthermore, since
$3\pi_2(v_2w)=L(2)\pi_2(v_0w)=20a\pi_2(v_2w)+12b\pi_2(v_w)$,
we have $20a+12b=3$. Solving the simultaneous equations, we have
$a=\frac{6}{41}$ and $b=\frac{1}{164}$.

For $\pi_3(v_{-2}u)$, there are $a,b\in \bC$ such that $\pi_3(v_{-2}u)=aL(-2)\pi_3(v_0u)+bL(-1)^2\pi_3(v_0u)$. Since $2L(-1)\pi_3(v_0u)=4\pi_3(v_{-1}u)=L(1)\pi_3(v_{-2}u)=3aL(-1)\pi_3(v_0u)+14bL(-1)\pi_3(v_0u)$, we have $3a+14b=2$.
Furthemore, since $5\pi_3(v_0u)=L(2)\pi_3(v_{-2}u)=24a\pi_3(v_0u)+18b\pi_3(v_0u)$,
we have $24a+18b=5$. Solving the simultaneous equations,
we have $a=\frac{17}{141}$ and $b=\frac{11}{94}$.
\prend

\begin{rmk}
The denominators are products of prime divisors of the order of the monster group. Tuit initially observed this phenomenon in \cite{T}.
\end{rmk}

\begin{rmk}
We do not use the following results for our proof but just for reference.
Set $X=\pi_{4}(v_0w)$ for $v\in V_2^p, w\in V_3^p$.
Then 
$\pi_{4} (v_{-2}w)=\frac{1}{8}L(-2)X+\frac{1}{16}L(-1)^2X$ and \\
$\pi_{4}(v_{-3}w)=\frac{1}{2^5\cdot 71}\{145L(-3)X+105L(-2)L(-1)X+\frac{79}{6}L(-1)^3X$. 
In particular, $L(-2)X\in L(-1)V_3+\pi_{4}((v\otimes e^{(1,1)})_0(u\otimes
\delta^{(1,-2)}(-1)e^{(1,2)}))$.
\end{rmk}

Conversely, using the above data, we can define physical states as follows:

\begin{defn}
For $0\leq k\leq 7$ and $Z\in V_k^p$, define
$$\begin{array}{rl}
\Phi^1_k(Z)=&\sum_{t=0}^{5-k}Q^{2,2}_k(-t)Z \otimes o_t(e^{(1,1)})e^{(2,2)}, \cr
\Phi^2_k(Z)=&\sum_{t=0}^{7-k}Q^{2,2}_k(-t)Z\otimes o_{t-1}(e^{(1,1)})\delta^{(1,-2)}_{(-1)}e^{(2,3)}, \qquad \mbox{ and}\cr
\Phi^3_k(Z)=&\sum_{t=0}^{7-k}Q^{2,3}_k(-t)Z\otimes o_t(e^{(1,1)})e^{(2,3)}.
\end{array}$$
\end{defn}

From the definition of $Q^{a,b}_k$, $\Phi^1_k$, $\Phi^2_k$, and $\Phi^3_k$, we have:

\begin{prop}
$\Phi^1(v,u)=\sum_{k=0}^{\infty} \Phi^1_k(\pi_k(v_{3-k}u))\in P^1_{(2,2)}, \quad
\Phi^2(v,u)=\sum_{k=0}^{\infty} \Phi^2_k(\pi_k(v_{3-k}u))\in P^1_{(2,3)}$, and
$\Phi^3(v,w)=\sum_{k=0}^{\infty}\Phi^3_k(\pi_k(v_{4-k}w))\in P^1_{(2,3)}$.
\end{prop}

Namely, $(v\otimes e^{(1,1)})_0(u\otimes e^{(1,1)}),
(v\otimes e^{(1,1)})_0(u\otimes \delta^{(1,-2)}_{(-1)}e^{(1,2)}),
(v\otimes e^{(1,1)})_0(w\otimes e^{(1,2)})$
are uniquely determined by
$$\begin{array}{ll}
(\pi_0(v_3u), \pi_2(v_1u), \pi_3(v_0u),\pi_4(v_{-1}u), \pi_5(v_{-2}u))&\hspace{-3mm}\in V_0\oplus (\oplus_{j=2}^5 V_j^p), \cr
(\pi_0(v_3u), \pi_2(v_1u), \pi_3(v_0u),\pi_4(v_{-1}u), \pi_5(v_{-2}u), \pi_6(v_{-3}u), \pi_7(v_{-4}u))
&\hspace{-3mm}\in V_0\oplus (\oplus_{j=2}^7V_j^p), \mbox{ and} \cr
(\pi_0(v_4w), \pi_2(v_2w), \pi_3(v_1w),\pi_4(v_0w),
\pi_5(v_{-1}w),\pi_6(v_{-2}w),\pi_7(v_{-3}w))
&\hspace{-3mm} \in V_0\oplus (\oplus_{j=2}^7V_j^p) ,
\end{array}$$
respectively. In particular, we have:

\begin{prop}\label{Once} \mbox{}\quad \\
(1) $(V_2^p\otimes e^{(1,1)})_0(V_2^p\otimes e^{(1,1)})$ covers $V_k^p$
modulo $\Rad(P^1)$ at most once for $k=0,2,3,4,5$. \\
(2) $(V_2^p\otimes e^{(1,1)})_0(V_2^p\otimes \delta^{(1,-2)}_{(-1)}e^{(1,2)})$
covers $V_k^p$ modulo $\Rad(P^1)$ at most once for $k=0,2,3,4,5,6,7$, and \\
(3) $(V_2^p\otimes e^{(1,1)})_0(V_3^p\otimes e^{(1,2)})$
covers $V_k^p$ modulo $\Rad(P^1)$ at most once for $k=0,2,3,4,5,6,7$.
\end{prop}

We can now prove the remaining parts of the proof of Lemma \ref{Afewcase}.
By the above lemma, the maximal bounds of multiplicities of
$V_0, V_2^p, V_3^p, V_4^p$, and $V_5^p$ in $B(V)^{(2,2)}$ and $B(V)^{(2,3)}$
are given. We also have $\dim B(V)^{(m,n)}=\dim V_{mn+1}$. Then, it is easy to check
$B(V)^{(2,2)}\cong V_{5}\cong B(V)^{(1,4)}$ and $B(V)^{(2,3)}\cong V_6$ as $G$-modules.
This completes the proof of Lemma \ref{Afewcase}. \\

Since $V_k^p$ are irreducible $G$-modules, we have:

\begin{lmm}
If there is $0\not=\alpha\in V_k^p$ such that $\Phi^j(\alpha)\in \Rad(P^1)$, then
$\Phi^j(V_k^p)\subseteq \Rad(P^1)$.
\end{lmm}

Therefore, Proposition \ref{Once} implies the following result.

\begin{prop}\label{Nointer}
Assume that $\{\Phi^1(v,u):v,u\in V_2^p\}$ covers $V_k^p$ modulo $\Rad(P^1)$.
If there are $v,u\in V_2^p$ such that $\pi_k(\Phi^1(v,u))\in \Rad(P^1)$,
then $\pi_k(\Phi^1(v,u))=0$, that is, $v_{3-k}u=0$.
Similar results hold for $\Phi^2$ and $\Phi^3$.
\end{prop}

Then we can translate a joint of Proposition \ref{Cover2} and \ref{Nointer} into the following:

\begin{lmm}\label{Cover3}
${\rm Span}_{\bC}\{(\pi_3(v_0u),\pi_5(v_{-2}u))\in V_3^p\oplus V_5^p\mid v,u\in V_2^p\}\cong V_5^p\oplus V_3^p$ and \\
${\rm Span}_{\bC}\{(\pi_0(v_3u), \pi_2(v_1u), \pi_2(v_2w), \pi_4(v_{-1}u), \pi_6(v_{-3}u), \pi_2(v'_2w), \pi_4(v'_0w), \pi_6(v_{-2}w))
\in V_0\oplus V_2^p\oplus V_4^p\oplus V_6^p\oplus V_2^p\oplus V_4^p\oplus V_6^p \mid v,v',u\in V_2^p, w\in V_3^p\}\cong V_0\oplus V_2^p\oplus V_2^p\oplus V_4^p\oplus V_4^p\oplus V_6^p$.
\end{lmm}

As a corollary of Lemma \ref{Direct} and \ref{Cover3}, we have:

\begin{lmm}\label{V5}
$V_5\subseteq C_2(V)$.
\end{lmm}

\pr
By Lemma \ref{Cover3}, for each $x\in V_3^p$, there are $v^i, u^i\in V_2^p$ such that
$0=\pi_5(\sum v^i_{-2}u^i)$ and $\pi_3(\sum v^i_0u^i)=x$.
Then since $\pi_3(\sum v^i_{-2}u^i)=\frac{17}{141}L(-2)x+\frac{11}{94}L(-1)^2x$ by Lemma \ref{Direct}, we have $\sum_i v^i_{-2}u^i-\frac{17}{141}L(-2)x\in
L(-1)V+\pi_0(V_5)+\pi_2(V_5)+\pi_4(V_5)\subseteq L(-1)V_4$,
which implies $L(-2)x\in C_2(V)$. By Proposition \ref{C2}, we have $V_5\subseteq C_2(V)$.
\prend

We will study $\pi_6(P^1_{(2,3)})$.

\begin{lmm}\label{V6}
For $Z\in V_6^p$, $\Phi^2(Z) \equiv 6\Phi^3(Z) \pmod{\Rad(P^1)}$.
\end{lmm}

\pr
For $Z\in V_6^p$, $Z\otimes e^{(2,3)}\in P^0_{(2,3)}$ and so
$L(-1)Z\otimes e^{(2,3)}+Z\otimes \delta^{(2,3)}_{(-1)}e^{(2,3)}\in \Rad(P^1_{(2,3)})$.
For $v,u\in V_2^p$ and $w\in V_3^p$, by (OP), we have
$$\begin{array}{l}
\pi_6((v\otimes e^{(1,1)})_0(w\otimes e^{(1,1)}))=\pi_6(v_{-2}w)\otimes \delta^{(1,1)}_{(-1)}e^{(2,3)}+\pi_6(v_{-3}w)\otimes e^{(2,3)} \cr
\mbox{}\quad=v_{-2}w\otimes \delta^{(1,1)}_{(-1)}e^{(2,3)}+\frac{5}{12}L(-1)v_{-2}w\otimes e^{(2,3)} \equiv \pi_6(v_{-2}w)\otimes (\delta^{(1,1)}_{(-1)}-\frac{5}{12}\delta^{(2,3)}_{(-1)})e^{(2,3)} \cr
\mbox{}\quad=\frac{1}{12}\pi_6(v_{-2}w)\otimes \delta^{(2,-3)}_{(-1)}e^{(2,3)} \pmod{\Rad(P^1)} \cr
\pi_6((v\otimes e^{(1,1)})_0(u\otimes \delta^{(1,-2)}_{(-1)}e^{(1,2)}))
=\pi_6(v_{-3}u)\otimes (\delta^{(1,-2)}_{(-1)}-\delta^{(1,1)}_{(-1)})e^{(2,3)}
-\pi_6(v_{-4}u)\otimes e^{(2,3)} \cr
\mbox{}\quad=\pi_6(v_{-3}u)\otimes \delta^{(0,-3)}_{(-1)}e^{(2,3)}-\frac{1}{2}L(-1)\pi_6(v_{-3}u)\otimes
e^{(2,3)}
\equiv \pi_6(v_{-3}u)\otimes (\delta^{(0,-3)}_{(-1)}+\frac{1}{2}\delta^{(2,3)}_{(-1)})e^{(2,3)}\cr
\mbox{}\quad=\frac{1}{2}\pi_6(v_{-3}u)\otimes \delta^{(2,-3)}_{(-1)}e^{(2,3)}
\pmod{\Rad(P^1)}
\end{array}$$
Therefore, we have the desired result.
\prend

\begin{lmm}
$L(-2)V_2^p+V_4^p \subseteq (V_2^p)_0(V_3)+L(-1)V_3$.
\end{lmm}

\pr
From Lemma \ref{Cover3}, we have:
$${\rm Span}_{\bC}\{(\pi_0(f_3u), \pi_2(f_1u), \pi_4(f_{-1}u), \pi_2(v_2w), \pi_4(v_0w))
\mid f,v,u\in V_2^p, w\in V_3^p\}$$
covers $V_0\oplus V_2^p\oplus V_4^p\oplus V_2^p\oplus V_4^p$ modulo $\Rad(P^1)$.
Since they covers $2$ copies of $V_2^p\oplus V_4^p$, Proposition \ref{Once} implies that
for any $(x,y)\in V_2^p\oplus V_4^p$, there are
$f^i, u^i, v^j\in V_2^p$ and $w^j\in V_3^p$ such that
$\pi_2(\sum_j v^j_2w^j)=x$, $\pi_4(\sum v^j_0w^j)=y$, and
$0=\pi_0(\sum f^i_3u^i)=\pi_2(\sum f^i_1u^i)=\pi_4(\sum f^i_{-1}u^i)$.
Then $\sum_j v^j_0w^j-\frac{3}{41}L(-2)x-\frac{21}{164}L(-1)^2x-y
\in \pi_3(V_4)\subseteq L(-1)V_3$, which implies
$L(-2)x+y\in L(-1)V_3+(V_2)_0(V_3)\subseteq (V_2)_0(V_3)$.
\prend

\begin{lmm}\label{V4}
$L(-2)^2V_2^p+L(-2)V_4^p\subseteq C_2(V)$.
\end{lmm}

\pr For each $\beta\in L(-2)V_2^p+V_4^p$, there are $v^i\in V_2$ and $w, w^i\in V_3$ such that $\beta=\sum v^i_0w^i+L(-1)w$.
Then $L(-2)\beta=\sum v^i_0(L(-2)w^i)+3\sum v^i_{-2}w^i+L(-2)L(-1)w\in C_2(V)$.
\prend

\begin{lmm}\label{V0}
$L(-2)^3V_0\subseteq C_2(V)$.
\end{lmm}

\pr
From Lemma \ref{Cover3} and \ref{V6}, we have that \vspace{-2mm}
$$
{\rm Span}_{\bC}\{(\pi_0(f_3u), \pi_2(f_1u+v_2w), \pi_4(f_{-1}u+v_0w),
\frac{\pi_6(6f_{-3}u+v_{-2}w)}{12}
\mid f,v,u\in V_2^p, w\in V_3^p\}$$
covers $V_0$.
Since they covers $2$ copies of $V_2^p\oplus V_4^p$, Proposition \ref{Once} implies
that for any $(x,y)\in V_2^p\oplus V_4^p$, there are
$f^i, u^i, v^j\in V_2^p$ and $w^j\in V_3^p$ such that
$\pi_0(\sum_j v^j_3u^j)=\1$ and
$0=\pi_2(\sum_j v^j_2w^j)=\pi_4(\sum v^j_0w^j)=\pi_2(\sum f^i_1u^i)=\pi_4(\sum f^i_{-1}u^i)=\pi_6(6\sum_j v^j_{-2}w^j+\sum f^i_{-3}u^i)$.
Namely, by Lemma \ref{Direct}, \vspace{-2mm}
$$\begin{array}{l}
196883(\sum_jv^j_{-2}w^j+\sum_i f^i_{-3}u^i)-3492L(-6)\1
-\frac{15623}{12}L(-4)L(-2)\1-\frac{1271}{2}L(-3)^2\1
\cr
\mbox{}\qquad -124L(-2)^3\1\in \pi_3(V_6)+\pi_5(V_6)\subseteq C_2(V).
\end{array}\vspace{-2mm}$$
We hence have $L(-2)^3\1\in C_2(V)$.
\prend

This completes the proof of Theorem 1.

\section{Appendix}
\begin{prop}
Let $V$ be a VOA satisfying the conditions of the uniqueness conjecture.
If there is a Poisson algebra isomorphism $\phi:V/C_2(V)\to V^{\natural}/C_2(V^{\natural})$ preserving the invariant bilinear form $\langle , \rangle$,
then $V\cong V^{\natural}$.
\end{prop}

\pr
Since $C_2(V)$ is a direct sum of homogeneous subspaces,
we may view $V/C_2(V)\cong \oplus_{m=0}^{\infty} (C_2(V)_m)^c$ 
as inner product spaces, and we have an isometry 
$$\tilde{\phi}:(\oplus_m C_2(V)_m)^c\to (\oplus_m C_2(V^{\natural})_m)^c,$$ 
where
$(C_2(V)_m)^c=\{v\in V_m \mid \langle v, u\rangle=0 \, 
\forall u\in C_2(V)\cap V_m\}$.
Let $\{\omega/\sqrt{12}, e^1,...,e^{196883}\}$ be an orthonormal basis of
$V_2$.
We note $\langle \omega, e^i_1e^j\rangle=\langle e^i_1\omega, e^j\rangle=\langle 2e^i,e^j\rangle=2\delta_{i,j}$.
For $i,j=1,...,196883$, since
$$\langle L(-2)e^s, e^i_{-1}e^j\rangle=\langle e^s, L(2)e^i_{-1}e^j\rangle
=\langle e^s, 4e^i_{1}e^j\rangle,$$
we have
$$\begin{array}{rl}
\langle \phi(e^s), \phi(e^i)_{1}\phi(e^j)\rangle
=&\frac{1}{4}\langle L(-2)\phi(e^s), \phi(e^i)_{-1}\phi(e^j)\rangle=\frac{1}{4}\langle L(-2)e^s, e^i_{-1}e^j\rangle \cr
=&\langle e^s, e^i_{1}e^j\rangle
\end{array}$$
for all $s,i,j$. 
Therefore, $\tilde{\phi}$ sends the $1$-product of $V_2$ to that of $V_2^{\natural}$.
In other words, $V_2$ is isomorphic to $V_2^{\natural}$ as Griess algebras.
Since their characters ${\rm ch}_V(\tau)$ and ${\rm ch}_{V^{\natural}}(\tau)$ are the same, we can obtain $V\cong V^{\natural}$ from \cite{DGL}.
\prend

\end{document}